\documentclass[10pt]{article}
\usepackage{amscd,amsfonts,amsmath,amssymb,amstext,amsthm}

\title{Injectivity of the Double Fibration Transform \\
for Cycle Spaces of Flag Domains}

\author{Alan T. Huckleberry\thanks{Research partially supported by the DFG
Schwerpunkt ``Global Methods in Complex Geometry''.}
\,\, \&\,\, Joseph A. Wolf
\thanks{Research partially supported by DFG Schwerpunkt ``Global Methods 
in Complex Geometry'' and by NSF Grant DMS 99--88643.\endgraf
\endgraf
}
}

\date{18 August 2003}
\newcommand{\s} {{\smallskip \noindent}}
\newcommand{\m} {{\medskip \noindent}}

\theoremstyle{plain}
\newtheorem{thm} [equation]{Theorem}
\newtheorem{lem} [equation]{Lemma}
\newtheorem{prop}[equation]{Proposition}
\newtheorem{cor} [equation]{Corollary}
\newtheorem{remk} [equation]{Remark}

\theoremstyle{definition}
\newtheorem*{defn} {Definition}

\theoremstyle{remark}

\theoremstyle {plain}

\begin{document}

\maketitle

\begin{abstract}
The basic setup consists of a complex flag manifold $Z=G/Q$ where $G$ is a
complex semisimple Lie group and $Q$ is a parabolic subgroup, an open orbit
$D = G_0(z) \subset Z$ where $G_0$ is a real form of $G$, and a
$G_0$--homogeneous holomorphic vector bundle $\mathbb E \to D$.
The topic here is the double fibration transform
${\cal P}: H^q(D;{\cal O}(\mathbb E)) \to H^0({\cal M}_D;{\cal O}(\mathbb E'))$
where $q$ is given by the geometry of $D$, ${\cal M}_D$ is the
cycle space of $D$, and $\mathbb E' \to {\cal M}_D$ is a certain naturally
derived holomorphic vector bundle.  Schubert intersection theory is used to
show that ${\cal P}$ is injective whenever $\mathbb E$ is sufficiently negative.\end{abstract}

\tableofcontents

\section {Introduction}
\setcounter{equation}{0}

Let $G_0$ be a non--compact real form of a complex semisimple Lie
group $G$ and consider its action on a compact $G$--homogeneous
projective algebraic manifold $Z=G/Q$.  It is of interest to
understand the $G_0$--representation theory associated to each
of its (finitely many) orbits.  In this note we restrict to
the case of an open orbit $D=G_0(z_0)$, and we often refer to such
open orbits as flag domains.

\m
The simplest example of this situation, where $G_0=SU(1,1)$,
$G=SL_2(\mathbb C)$ and $Z=\mathbb P_1(\mathbb C)$, is at first
sight perhaps somewhat misleading.  In this case $D$ is either the
set of of negative or positive lines and therefore is biholomorphic to the
unit disk $\Delta =\{z\in \mathbb C:\vert z\vert <1\}$.  The naturally
associated $G_0$--representations can be regarded as being in
$L^2$--spaces of holomorphic functions on ${\cal O}(\Delta )$, or their
complex conjugates.

\m
Unless $G_0$ is of Hermitian type and $D$ is biholomorphically equivalent
to the associated bounded symmetric domain as in the above example, the
holomorphic functions do not separate the points of $D$; in fact in
most cases ${\cal O}(D)\cong \mathbb C$.  One explanation for this is
that $D$ generally contains positive--dimensional compact analytic subsets
which are in fact closely related to the representation theory at hand.
These arise initially as orbits  of maximal compact subgroups $K_0$
of $G_0$: {\it There is a unique complex $K_0$--orbit $C_0$ in $D$}
\cite {W1}.

\m
Of course $C_0$ can be just a point, which is exactly the case when
$D$ is a bounded Hermitian symmetric domain, but in general it
has positive dimension.  

\s
For example, if $G_0=SL_3(\mathbb R)$,
$G=SL_3(\mathbb C)$ and $Z=\mathbb P_2(\mathbb C)$, then the closed
$G_0$--orbit in $Z$ is the set of real points 
$Z_\mathbb R=\mathbb P_2(\mathbb R)$ and its complement 
$D=Z\setminus Z_\mathbb R$  is the only other orbit.  If $K_0$ is
chosen to be the real orthogonal group $K_0=SO_3(\mathbb R)$, then
\begin {equation*}
C_0=\{ [z_0:z_1:z_2]:z_0^2+z_1^2+z_2^2=0\}
\end {equation*}
is the standard quadric.  Here one easily checks that 
${\cal O}(D)\cong \mathbb C$ and, in view of basic results of Andreotti 
and Grauert \cite {AnG},
looks for $G_0$--representations in Dolbeault cohomology
$H^1(D,\mathbb E)$, where $\mathbb E \to D$ is a sufficiently negative 
holomorphic vector bundle.
  
\m
Cohomology classes, e.g., classes of bundle valued differential forms,
are technically and psychologically more difficult to handle than holomorphic
functions or sections of vector bundles.  Thus, for $q:=\dim_\mathbb CC_0$
one is led to consider the space ${\cal C}^q(D)$ of $q$--dimensional
cycles in $D$, where $C_0$ can be considered as a point.

\s
An element $C\in {\cal C}^q(D)$ is a linear combination,
$C=n_1C_1+\ldots n_mC_m$, with $n_j\in \mathbb N^{>0}$ and where
$C_j$ is a $q$--dimensional irreducible compact analytic subset of $D$.
It is also necessary to consider ${\cal C}^q(Z)$, where 
the $C_j$ are not required to be contained in $D$.

\m
The cycle space ${\cal C}^q(X)$ of a complex space $X$ has a
canonical structure of (locally) finite--dimensional complex
space (see \cite{B} and, e.g., \cite {GPR}, for this and other
basic properties).  In the case at hand, where $Z$ is projective
algebraic, ${\cal C}^q(Z)$ is often referred to as the Chow variety; in
particular, its irreducible components are projective algebraic
varieties.

\s
To simplify the notation we redefine ${\cal C}^q(Z)$ to be the
topological component containing the base cycle $C_0$\,.  It
contains ${\cal C}^q(D)$ as an open semialgebraic subset.  This 
is well defined independent of the choice of the maximal compact 
subgroup $K_0$.

\m
It is known that the induced action of $G$ on ${\cal C}^q(Z)$ is
algebraic (see, e.g., \cite {Hei} for a detailed proof) and therefore
the orbit ${\cal M}_Z :=G.C_0$ is Zariski open in its closure $X$ in
${\cal C}^q(Z)$.  It should be noted that ${\cal M}_Z $ is a spherical
homogeneous space of the reductive group $G$ and therefore the
Luna--Vust theory of compactifications \cite {BLV} applies.
It would be extremely interesting to determine the $G$--varieties
$X$ which occur in this way.

\m
In general it may be difficult to understand the full cycle space
${\cal C}^q(D)$ and therefore one cuts down to a simpler space which
is more closely related to the group actions at hand.  That simpler
space is

\begin {equation*}
{\cal M}_D: \text{ topological component of } C_0 \text{ in }
{\cal M}_Z \cap {\cal C}^q(D).
\end {equation*}

We have the incidence space

\begin {equation*}
{\mathfrak X_D} =
\{ (z,C)\in D\times {\cal M}_D:z\in C\}.
\end {equation*}
and the projections
$\mu:{\mathfrak X_D}\to D$ by $(z,C)\mapsto z$, and
$\nu:{\mathfrak X_D}\to {\cal M}_D$ by $(z,C)\mapsto C$.

\m
If $\mathbb E \to D$ is a holomorphic vector bundle, we may lift it
to $\mu^*\mathbb E \to \mathfrak X_D$ consider the associated $\nu $--direct
image sheaves on ${\cal M}_D$\,.  In this way Dolbeault cohomology spaces 
$H^q(D, \mathbb E)$ are transformed to the level of sections of holomorphic
vector bundles $\mathbb E' \to {\cal M}_D$\,.  This {\it double fibration 
transform} is explained in
detail in the next section.  Since ${\cal M}_D$ is now known to
be a Stein domain (\cite {W2}, \cite {HW}), our somewhat technical initial
setting is transformed to one that is more tractable. 

\m
In recent work we developed complex geometric methods aimed at 
describing the cycle spaces ${\cal M}_D$.  For example, for a
fixed real form, the space ${\cal M}_D$ is essentially always
biholomorphically equivalent to a universal domain ${\cal U}$
which is defined independent of $D$ and $Z$ (\cite {HW},\cite {FH}).

\s
On the other hand, here we prove that ${\cal M}_D$ possesses canonically
defined holomorphic fibrations which do indeed depend on $D$ and $Z$ and
which give it interesting refined structure.  As a consequence we
in particular show that the double fibration transform is injective,
a fact that should be useful for representation theoretic considerations.

\m
The present paper is organized as follows. In $\S 2$ we recall the basics
of the double fibration transformation.  The recently proved results
of (\cite {HW}) and (\cite {FH}) are summarized in $\S 3$ and in
$\S 4$ we give details on the method of Schubert slices. In the final section
these slices are used to define the fibrations of 
${\cal M}_D$ which were mentioned above.  The injectivity of the
double fibration transform is a consequence of the topological triviality
of these bundles. 

\section {Basics of the Double Fibration Transform}
\setcounter{equation}{0}

Let $D$ be a complex manifold (later it will be an open orbit of
a real reductive group $G_0$ on a complex flag manifold $Z = G/Q$
of its complexification ).  We suppose that $D$ fits into what we loosely
call a {\sl holomorphic double fibration}.  This means that there are complex
manifolds ${\cal M}$ and $\mathfrak X$ with maps
\begin{gather} \label{gen_doublefibration}
\setlength{\unitlength}{.08 cm}
\begin{picture}(180,18)
\put(81,15){$\mathfrak X$}
\put(74,12){$\mu$}
\put(69,1){$D$}
\put(90,12){$\nu$}
\put(95,1){$\cal M$}
\put(80,13){\vector(-1,-1){6}}
\put(87,13){\vector(1,-1){6}}
\end{picture}
\end{gather}
where $\mu$ is a holomorphic submersion and $\nu$ is a proper holomorphic
map which is a locally trivial bundle.  Given a locally free coherent 
analytic sheaf ${\cal E} \to D$ we construct a locally free coherent analytic 
sheaf ${\cal E'} \to {\cal M}$ and a transform
\begin{equation}
{\mathcal P} : H^q(D;{\cal E}) \to H^0({\cal M};{\cal E'})
\end{equation}
under mild conditions on (\ref{gen_doublefibration}). 
This construction is fairly standard,
but we need several results specific to the case of flag domains.
\m

\noindent {\bf Pull--back.}
\m

The first step is to pull cohomology back from $D$ to $\mathfrak X$.  Let
$\mu^{-1}({\cal E}) \to \mathfrak X$ denote the inverse image sheaf.
For every integer $r \geqq 0$ there is a natural map
\begin{equation} \label{pull1}
\mu^{(r)} : H^r(D;{\cal E}) \to H^r(\mathfrak X; \mu^{-1}({\cal E}))
\end{equation}
given on the \v Cech cocycle level by $\mu^{(r)}(c)(\sigma) = c(\mu(\sigma))$
where $c \in Z^r(D;\cal E)$ and where
$\sigma = (w_0, \dots , w_r)$ is a simplex.  For $q \geqq 0$
we consider the Buchdahl $q$--condition on the fiber $F$ of 
$\mu: \mathfrak X \to D$:
\begin{equation} \label{buchdahl_conditions}
 F \text{ is connected and }
H^r(F;\mathbb C) = 0 \text{ for } 1 \leqq r \leqq q-1.
\end{equation}

\begin{prop} \label{buchdahl_theorem} {\rm (See \cite{Bu}.)}
Fix $q \geqq 0$.
If {\em (\ref{buchdahl_conditions})} holds, then
{\rm (\ref{pull1})} is an isomorphism
for $r \leqq q-1$ and is injective for $r = q$.
If the fibers of $\mu$ are cohomologically acyclic then
{\rm (\ref{pull1})} is an isomorphism for all $r$.
\end{prop}

As usual, ${\cal O} _X \to X$ denotes the structure sheaf of a complex
manifold $X$ and ${\cal O}(\mathbb E) \to X$ denotes the sheaf of germs of
holomorphic sections of a holomorphic vector bundle $\mathbb E \to X$.
Let $\mu^*({\cal E}) := 
\mu^{-1}({\cal E}) \widehat{\otimes}_{\mu^{-1}({\cal O}_D)}
{\cal O}_{\mathfrak X} \to \mathfrak X$ denote the pull--back sheaf.  
It is a coherent analytic sheaf of ${\cal O}_{\mathfrak X}$--modules.  
Now $[\sigma] \mapsto [\sigma] \otimes 1$ defines a map
$i : \mu^{-1}({\cal E}) \to \mu^*({\cal E})$
which in turn specifies maps in cohomology, the coefficient morphisms
\begin{equation}\label{pull2}
i_p : H^p(\mathfrak X; \mu^{-1}({\cal E})) \to H^p(\mathfrak X; \mu^*({\cal E})) \ \
\text{ for } p \geqq 0.
\end{equation}
Our natural pull--back maps are the compositions
$j^{(p)} = i_p \cdot \mu^{(p)}$ of (\ref{pull1}) and (\ref{pull2}):
\begin{equation} \label{pull}
j^{(p)} : H^p(D;{\cal E}) \to H^p(\mathfrak X;\mu^*({\cal E})) \ \ \text{ for } p \geqq 0.
\end{equation}
We have $ {\cal E} = {\cal O}(\mathbb E)$ for some holomorphic vector bundle
$\mathbb E \to D$, because we assumed $ {\cal E} \to D$ locally free.
Thus $\mu^*({\cal E}) = {\cal O}(\mu^*(\mathbb E))$ and we
realize these sheaf cohomologies as Dolbeault cohomologies.  In the context
of Dolbeault cohomology, the pull--back maps (\ref{pull}) are given by 
pulling back $[\omega] \mapsto [\mu^*(\omega)]$ on the level of differential 
forms.
\medskip

\noindent {\bf Push--down.}
\medskip

In order to push the $H^q(\mathfrak X;\mu^*({\cal E}))$ down to
${\cal M}$ we assume that
\begin{equation} \label{proper_stein}
{\cal M} \text{ is a Stein manifold}.
\end{equation}
Since $\nu : {\mathfrak X} \to {\cal M}$ was assumed proper, we have the Leray 
direct image sheaves ${\cal R}^p(\mu^*({\cal E})) \to {\cal M}$.  Those
sheaves are coherent \cite{GR}.  As ${\cal M}$ is Stein
\begin{equation}
H^q({\cal M}; {\cal R}^p({\cal E})) = 0 \ \ \text{ for } \ \ p \geqq 0 \text{ and } q > 0.
\end{equation}
Thus the Leray spectral sequence of $\nu : {\mathfrak X} \to {\cal M}$ 
collapses and gives
\begin{equation} \label{collapse}
H^p(\mathfrak X; \mu^*({\cal E})) \cong H^0({\cal M};{\cal R}^p(\mu^*({\cal E}))).
\end{equation}

\begin{defn} 
The {\em double fibration transform} for the double
fibration {\rm (\ref{gen_doublefibration})} is the composition
\begin{equation} \label{transform}
{\mathcal P}: H^p(D; {\cal E}) \to H^0({\cal M}; {\cal R}^p(\mu^*({\cal E})))
\end{equation}
of the maps {\rm (\ref{pull})} and {\rm (\ref{collapse})}.

\end{defn}

In order that the double fibration transform (\ref{transform})
be useful, one wants two conditions to be satisfied.  They are
\begin{gather}
\phantom{X} {\mathcal P}: H^p(D; {\cal E}) \to H^0({\cal M}; {\cal R}^p(\mu^*({\cal E})))
\text{ should be injective, and}
\label{want_inj}\\
\text{there should be an explicit description of the image of }
{\mathcal P} \label{want_image}.
\end{gather}
Assuming (\ref{proper_stein}),
injectivity of ${\mathcal P}$ is equivalent to injectivity of $j^{(p)}$ in
(\ref{pull}).  The most general way to approach this is the combination of
vanishing and negativity in Theorem \ref{general_injective} below, taking the
Buchdahl conditions (\ref{buchdahl_conditions}) into consideration.
\m

Given the setting of (\ref{proper_stein}) our attack on the 
injectivity question uses a
spectral sequence argument for the relative de Rham complex
of the holomorphic submersion $\mu:\mathfrak X \rightarrow D$.  
See \cite{WZ} for the details.  The end result is

\begin{thm} \label{general_injective} 
Let ${\cal E}={\cal O}(\mathbb E)$ for some
holomorphic vector bundle $\mathbb E \to D$.
Fix $q \geqq 0$.  Suppose that
the fiber $F$ of $\mu : \mathfrak X \to D$ is connected and satisfies
{\rm (\ref{buchdahl_conditions})}.  Assume {\rm (\ref{proper_stein})}
that ${\cal M}$ is Stein. 
Suppose that $H^p(C; \Omega_\mu^r(\mathbb E)|_C) = 0$ for
$p < q$, and $r \geqq 1$ for every fiber $C$ of $\nu : {\mathfrak X} \to {\cal M}$
where $\Omega_\mu^r(\mathbb E) \to \mathfrak X$ denotes the sheaf of relative
$\mu^*\mathbb E$--valued holomorphic $r$--forms on $\mathfrak X$ with respect 
to $\mu : \mathfrak X \to D$.  Then
${\mathcal P} :H^q(D;{\cal E}) \to H^0({\cal M}; {\cal R}^q(\mu^*{\cal E}))$ is injective.
\end{thm}
\medskip

\noindent {\bf Flag domain case.}
\medskip

In the cases of interest to us, $D$ will be a flag domain, we will 
have ${\cal E}={\cal O}(\mathbb E)$ 
as in Theorem \ref{general_injective}, and the transform ${\mathcal P}$ will have an
explicit formula.  The  Leray derived sheaf will be given by
\begin{equation} \label{derived_bundle_gen}
\begin{aligned} {\cal R}^q&(\mu^*({\cal O}(\mathbb E))) 
	= {\cal O}(\mathbb E') \text{ where } \\
& \mathbb E' \to {\cal M} \text{ has fiber }
H^q(\nu^{-1}(C); {\cal O}(\mu^*(\mathbb E)|_{\nu^{-1}(C)})) \text{ at } C.
\end{aligned}
\end{equation}
Then ${\mathcal P}$ will be given on the level of Dolbeault cohomology, as follows.
Let $\omega$ be an $\mathbb E$--valued $(0,q)$--form on $D$  and 
$[\omega] \in H_{\overline{\partial}}^q(D,\mathbb E)$ its Dolbeault class.  
Then
\begin{equation*}
\begin{aligned}
{\mathcal P}([\omega]) &\text{ is the holomorphic section of }
\mathbb E' \to {\cal M} \\
&\text{ whose value } {\mathcal P}([\omega])(C) \text{ at }
C \in {\cal M} \text{ is } [\mu^*(\omega)|_{\nu^{-1}(C)}].
\end{aligned}
\end{equation*}
In other words,
\begin{equation}
{\mathcal P}([\omega])(C) = [\mu^*(\omega)|_{\nu^{-1}(C)}] \in
H^0_{\overline{\partial}} ({\cal M}; \mathbb E')
\end{equation}
This is most conveniently interpreted by viewing ${\mathcal P}([\omega])(C)$ as
the Dolbeault class of $\omega|_{C}$, and by viewing $C \mapsto
[\omega|_{C}]$ as a holomorphic section of the holomorphic vector
bundle $\mathbb E' \to {\cal M}$.
\medskip

Now let $D = G_0(z_0)$ be an open orbit in the complex flag manifold
$Z = G/Q$, and ${\cal M}$ is replaced by the cycle space ${\cal M}_D$.  
Our double fibration (\ref{gen_doublefibration}) is replaced by
\begin{gather} \label{flag_doublefibration}
\setlength{\unitlength}{.08 cm}
\begin{picture}(120,18)
\put(61,15){$\mathfrak X_D$}
\put(50,12){$\mu$}
\put(46,1){$D$}
\put(75,12){$\nu$}
\put(77,1){${\cal M}_D$}
\put(58,13){\vector(-1,-1){6}}
\put(69,13){\vector(1,-1){6}}
\end{picture}
\end{gather}
where $\mathfrak X_D := \{(z,C) \in D \times {\cal M}_D \mid z \in C \}$ is the
incidence space.  Given a $G_0$--homogeneous holomorphic vector bundle
$\mathbb E \to D$, and the number $q = \dim_\mathbb C C_0$\,, we will see in
Section \ref{der-sheaves} that the Leray derived sheaf
involved in the double fibration transform satisfies 
(\ref{derived_bundle_gen}).   Here (\ref{derived_bundle_gen}) will become a
little bit more explicit and take the form
\begin{equation}\label{derived_bundle}
\begin{aligned}
{\cal R}^q(\mu^*({\cal O}(\mathbb E))) &= {\cal O}(\mathbb E') 
	\text{ where } \\
&\mathbb E' \to {\cal M}_D \text{ has fiber } 
	H^q(C; {\cal O}(\mathbb E|_{C})) \text{ at } C \in {\cal M}_D\ .
\end{aligned}
\end{equation}
Evidently, $\mathbb E' \to {\cal M}_D$ will be globally $G_0$--homogeneous. 
It cannot be $G$--homogeneous unless ${\cal M}_D$ is $G$--invariant, and that
only happens in the degenerate case where $G_0$ is transitive on $Z$.  
However, in Section \ref{der-sheaves} we will see that 
$\mathbb E' \to {\cal M}_D$ is the restriction of a $G$--homogenous
holomorphic vector bundle $\widetilde{\mathbb E}' \to {\cal M_Z}$\,,
and in particular (\ref{derived_bundle_gen}) is satisfied.  In any  
case, $H^q(C; {\cal O}(\mathbb E|_{C}))$ can be
calculated from the Bott--Borel--Weil Theorem.
Thus ${\cal R}^q(\mu^*({\cal O}(\mathbb E)))$ will be
given explicitly by (\ref{derived_bundle}) in the flag domain case.
\medskip

Using methods of complex geometry as described in Section \ref{cycle-space} 
below it was 
shown (\cite{HW}  plus \cite{FH}) that ${\cal M}_D$ is biholomorphically
equivalent to a certain universal domain ${\cal U}$ or to a bounded
symmetric domain.  It follows in general that ${\cal M}_D$ is a contractible 
Stein manifold, so $\mathbb E' \to {\cal M}_D$ is holomorphically
trivial.   We will use those same methods in Section \ref{schubert} to 
construct 
certain holomorphic fibrations of the ${\cal M}_D$ and use 
those ``Schubert fibrations'' to show that $F$ satisfies 
(\ref{buchdahl_conditions}) for all $q$.  That is how we will prove that 
the double fibration transforms are injective.  Thus, in the flag
domain case, we will have a complete answer to (\ref{want_inj}) and some
progress toward (\ref{want_image}).

\section {A Computable Description of ${\cal M}_D$} \label{cycle-space}
\setcounter{equation}{0}

In order to understand the structure of $Z$, $D$ and ${\cal M}_D$ we may
assume that $G_0$ is simple, because $G_0$ is local direct product of simple
groups, and $Z$, $D$ and ${\cal M}_D$ break up as global direct products
along the local direct product decomposition of $G_0$\,.  From this point
on $G_0$ is simple unless we say otherwise.  We also assume that $G_0$
is not compact, and we avoid the two trivial noncompact cases (see \cite{W4}),
where $G_0$ acts transitively on $Z$ so that ${\cal M}_D$ is reduced 
to a single point. 
\m

As above let ${\cal M}_Z :=G.C_0$. Since $C_0$ is a complex manifold,
the isotropy group $G_{C_0}$ contains the complexification $K$ of $K_0$\,,
which is a closed complex subgroup of $G$. If $G_0$ is not Hermitian, then
$\mathfrak k$ is a maximal subalgebra of $\mathfrak g$, so $G_{C_0}/K$ is finite.  
In this nonhermitian case, without further discussion,
we replace ${\cal M}_Z $ by the finite cover $G/K$.  There is no loss of 
generality in doing this because, as we will see later,
${\cal M}_D$ pulls back biholomorphically to a domain in $G/K$
under the finite covering $G/K \to G.C_0$\,. 

\m
If $G_0$ is Hermitian, then the symmetric space $G_0/K_0$ possesses
two invariant complex structures, ${\cal B}$ and $\overline {\cal B}$,
as bounded Hermitian symmetric domains.  These are realized
as open orbits $D$ in their compact duals $G/P$ and $G/\bar P$
respectively.  In these cases the base cycles are just the
$K_0$ fixed points and clearly $G_{C_0}$ is either $P$ or $\overline P$
in such a situation; in particular, as opposed to being the
affine homogeneous space $G/K$, the space ${\cal M}_Z $ is a compact
homogeneous manifold.

\m
The above mentioned phenomenon is characterized by $G_{C_0}$ being
either $P$ or $\bar P$. It occurs in more interesting situations
than just that where $D$ itself is a bounded domain, but it should
be regarded as an exceptional case which is completely understood
\cite {W2}. In particular, in any such example the cycle space
${\cal M}_D$ is either ${\cal B}$ or $\overline {\cal B}$.

\m
In all other cases it was recently shown that ${\cal M}_D$
is naturally biholomorphic to a universal domain
${\cal U}\subset {\cal M}_Z =G/K$ (\cite {FH}, \cite {HW}).  It should
be emphasized that ${\cal U}$, which we now define, depends
{\it only} on the real form $G_0$, and not on $D$.

\m
The domain ${\cal U}$ can be defined in a number of different ways.
We choose the historical starting point which is of a differential
geometric nature.

\m
Let $M$ denote the Riemannian symmetric space $G_0/K_0$ of negative
curvature and consider its tangent bundle $TM$.  As usual let
$\theta $ be a Cartan involution of ${\mathfrak g}$ that commutes
with complex conjugation over its real form ${\mathfrak g}_0$\,.  
$\theta$  defines the
Cartan decompositions ${\mathfrak g}={\mathfrak g_u}\oplus i{\mathfrak g_u}$
and ${\mathfrak g}_0={\mathfrak k}_0\oplus {\mathfrak s}_0$\,, where
$\mathfrak g_u$ is the compact real form of $\mathfrak g$ with
$\mathfrak k_0 = \mathfrak g_u \cap \mathfrak g_0$\,.

\m
We identify ${\mathfrak s}_0$ with $TM_{x_0}$, where $x_0$
is the base point with $G_0$--isotropy $K_0$, and regard $TM$
as the homogeneous bundle $G_0\times _{K_0}{\mathfrak s}_0$.
One defines the {\it polar coordinates} mapping by

\begin {equation*}
{\Pi}:TM\to G/K, \ ([g_0,\xi])\mapsto g_0 \exp(i\xi ).x_0.
\end {equation*}

\m
Clearly ${\Pi}$ is a diffeomorphism in a neighborhood of the
$0$--section and thus it is of interest to consider the canonically
defined domain

\begin {equation*}
\Omega _{max}:=\{ v\in TM:\text{rank} ({\Pi}_*(v))= \dim\, TM\}^0.
\end {equation*}

Here the connected component is that which contains the $0$--section.

\m
It turns out that $\Omega _{max}$ is determined by differential
geometric properties of the compact dual $N:=G_u.{x_0}=G_u/K_0$, 
where $G_u$ is the maximal compact subgroup of $G$ defined
by ${\mathfrak g_u}$.  For this it is essential that the geodesics
emanating from $x_0$ in $N$ are just orbits $ \exp(i\xi ).x_0$,
$\xi \in {\mathfrak s}_0$, of $1$--parameter groups.

\s
Let $\frac{1}{2}N$ denote the set of points in $N$ which are at most
halfway from $x_0$ to the cut point locus and define
$\Omega _C=G_0.\frac{1}{2}N$.

\begin {thm} {\rm (Crittenden \cite {C})}
The polar coordinates mapping ${\Pi}$ restricts to a diffeomorphism
from $\Omega _{max}$ to $\Omega _C$.
\end {thm}

The domain $\Omega _C$ can be computed in an elementary way.  For
this regard $K_0$ as acting on ${\mathfrak s}_0$ by the adjoint
representation, let ${\mathfrak a}_0$ be a maximal Abelian subalgebra
in ${\mathfrak s}_0$ and recall that 
$K_0.{\mathfrak a}_0={\mathfrak s}_0$.  Thus $\Omega _C$ is determined
by a domain in ${\mathfrak a}_0$.  This is computed as follows.

\s
Since ${\mathfrak a}_0$ acts on ${\mathfrak g}_0$ as a commutative algebra
of self--adjoint transformations, the associated joint eigenvalues 
$\alpha \in {\mathfrak a}_0^*$ are real valued.  (The nonzero ones are the
$(\mathfrak g_0, \mathfrak a_0)$--roots or {\em restricted} roots.)\,  Consider
the convex polygon

\begin {equation*}
V = \{ \xi \in {\mathfrak a}_0 \mid |\alpha(\xi)| < \tfrac{\pi}{2}
	\text{ for all restricted roots } \alpha \}
\end {equation*}

and define

\begin {equation*}
{\cal U}:=G_0. \exp(iV).x_0.
\end {equation*}

\begin {thm} {\rm (\cite {C}, \cite {AkG})}
$\Omega _C={\cal U}$.

\end {thm}

\begin {remk}
{\em The domain ${\cal U}$, which is indeed explicitly computable, was
first brought to our attention by the work in \cite {AkG}.  As
a consequence we originally denoted it by $\Omega _{AG}$.  It turns
out that it is naturally equivalent to a number of other domains, 
including the cycle spaces, which are defined from a variety of
viewpoints.  So now, unless we have a particular construction
in mind, we denote it by ${\cal U}$ to underline its universal
character.}
\end {remk}

\begin {prop} \label {contractible}
The domain ${\cal U}$ is contractible.
\end {prop}

\begin {proof}
For this we regard ${\cal U}$ as being contained
in $TM$.  Since $V$ is a (polyhedral) domain in ${\mathfrak s}_0$
star--shaped from $0$, scalar multiplication
$\varphi _t:TM\times TM$, $v\mapsto tv$, stabilizes $V$, and for
$0\le t\le 1$ defines a strong deformation retraction of 
${\cal U}$ to the $0$--section $M$.  That $0$--section is also
contractible.
\end {proof}

Complex geometric properties of ${\cal U}$ are also of importance.  
These include the fact that ${\cal U}$ is a Stein domain in 
${\cal M_Z}$ (\cite{BHH}; also see \cite{Ba}, \cite{H}, \cite{HW},
\cite{GK}, \cite{GM}). 
\medskip

For later reference let us state the main 
points in the context of cycle spaces.

\begin {thm} \label {cycle spaces}
If $D$ is not one of the above exceptions where 
${\cal M}_D$ is either a single point or ${\cal B}$ or $\bar {\cal B}$, then
${\cal M}_D={\cal U}$.  Thus in all cases ${\cal M}_D$ is contractible
and Stein.
\end {thm}

\begin{cor} \label {cycle-maps}
Let $P \subset Q$ be parabolic subgroups of $G$.  Let $\pi$ denote the
natural projection $1P \mapsto 1Q$ of $W = G/P$ onto $Z = G/Q$.  Suppose
that the flag domain $D \subset Z$ is not one of the above exceptions where
${\cal M}_D$ is either a single point or ${\cal B}$ or $\bar {\cal B}$.
Let $\widetilde{D} \subset W$ be a flag domain such that $\pi(\widetilde{D})
= D$.  Then $\pi$ induces a holomorphic diffeomorphism of 
${\cal M}_{\widetilde{D}}$ onto ${\cal M}_D$\,.
\end{cor}

\section{Globalization of the Bundles} 
\label{der-sheaves}
\setcounter{equation}{0}

We will show that various bundles are restrictions of bundles homogeneous 
under the complex group $G$, and use that to carry bundles over $D$ to
bundles over ${\cal M}_D$\,.  For that we need an old result on homogeneous
holomorphic vector bundles from \cite[Section 3]{TW}.
\m

$M = A_0/B_0$ be a homogeneous complex manifold.  Let $p: A_0 \to M$ denote
the natural projection.  View the Lie algebra ${\mathfrak a}_0$\,, and thus its
complexification ${\mathfrak a}$, as Lie algebras of holomorphic vector fields
on $M$.  Let $m \in M$ denote the base point $1B_0$ and define
${\mathfrak p} = \{\xi \in {\mathfrak a} \mid \xi_m = 0\}$.  Then
${\mathfrak p}$ is an ${\rm Ad}(B_0)$--stable complex subalgebra of 
${\mathfrak a}$ such that 
${\mathfrak a} = {\mathfrak p} + \overline{\mathfrak p}$ and
${\mathfrak b} = {\mathfrak p} \cap \overline{\mathfrak p}$.  Here
${\mathfrak b}$ is the complexification of the Lie algebra ${\mathfrak b}_0$
of $B_0$\,.
\s

Let $\chi$ be a continuous representation of $B_0$ on a finite dimensional
complex vector space $E = E_\chi$\,.  By {\em extension of} $\chi$ to 
${\mathfrak p}$ we mean a Lie algebra representation $\lambda$ of 
${\mathfrak p}$ on $E$ such that $\lambda|_{{\mathfrak b}_0} = d\chi$ and
$\lambda(\text{\rm Ad}(b)\xi) = \chi(b)\lambda(\xi)\chi(b)^{-1}$ for all
$b \in B_0$ and $\xi \in {\mathfrak p}$.  Thus an extension
of $\chi$ to ${\mathfrak p}$ is a $({\mathfrak p},B_0)$--module structure
on $E_\chi$\,.
\s

The representation $\chi$ defines a real analytic, $A_0$--homogeneous,
complex vector bundle ${\mathbb E}_\chi \to M =  A_0/B_0$\,, by
${\mathbb E}_\chi = A_0 \times_{B_0} E_\chi$\,.  We identify
local sections $s : U \to {\mathbb E}_\chi$ with functions
$f_s : p^{-1}(U) \to E_\chi$ such that $f_s(gb) = \chi(b)^{-1}f_s(b)$.
\s

The holomorphic vector bundle structures on ${\mathbb E}_\chi \to M$ 
are given as follows.

\begin{thm} \label{t-w}
{\rm \cite[Theorem 3.6]{TW}}  The structures of $A_0$--homogeneous
holomorphic vector bundle on ${\mathbb E}_\chi \to M$ are on one to one
correspondence with the extensions $\lambda$ of $\chi$ from $B_0$ to
${\mathfrak p}$.  The structure corresponding to $\lambda$ is the one for which
the holomorphic sections $s$ over any open set $U \subset M$ are characterized 
by $\xi\cdot f_s + \lambda(\xi)f_s = 0$ on $p^{-1}(U)$ for all 
$\xi \in {\mathfrak p}$.
\end{thm}

Now we return to the flag domain setting and extend the double
fibration (\ref{flag_doublefibration}) to
\begin{gather} \label{global_doublefibration}
\setlength{\unitlength}{.08 cm}
\begin{picture}(120,18)
\put(61,15){$\mathfrak X_Z$}
\put(50,12){$\widetilde{\mu}$}
\put(46,1){$Z$}
\put(75,12){$\widetilde{\nu}$}
\put(77,1){${\cal M}_Z$}
\put(58,13){\vector(-1,-1){6}}
\put(69,13){\vector(1,-1){6}}
\end{picture}
\end{gather}
where ${\cal M}_Z = \{gC_0 \mid g \in G\}$ and 
$\mathfrak X_Z := \{(z,C) \in Z \times {\cal M}_Z \mid z \in C \}$ is the
incidence space.  Evidently ${\cal M}_D$ is an open submanifold of
${\cal M}_Z$\,, and $\mu$ and $\nu$ are the
respective restrictions of $\widetilde{\mu}$ and $\widetilde{\nu}$.
\m

\begin{thm} \label{derived-bundle}  Let $D$ be an open
$G_0$--orbit on $Z$, let ${\mathbb E} \to D$ be a $G_0$--homogeneous 
holomorphic vector bundle, and let $q \geqq 0$.  Suppose $G_0 \subset G$.
Then {\rm (1)} 
${\mathbb E} \to D$ is the restriction of a $G$--homogeneous holomorphic
vector bundle $\widetilde{\mathbb E} \to Z$, 
{\rm (2)} the Leray derived sheaf for $\widetilde{\nu}$ is given by
${\cal R}^q({\cal O}(\widetilde{\mu}^*\widetilde{\mathbb E})) = 
{\cal O}(\widetilde{\mathbb E}')$ where
$\widetilde{\mathbb E}' \to {\cal M}_Z$ is the $G$--homogeneous,
holomorphic vector bundle with fiber $H^q(C;{\cal O}(\widetilde{\mathbb E}|_C))$ 
over $C \in {\cal M}_Z$\,, and {\rm (3)} the Leray derived sheaf for $\nu$
is given by ${\cal R}^q({\cal O}(\mu^*{\mathbb E})) =
{\cal O}({\mathbb E}')$ where ${\mathbb E}'$ is the restriction of
$\widetilde{\mathbb E}'$ to ${\cal M}_D$\,.
\end{thm}
 
\begin{proof}
We translate the result of Theorem \ref{t-w} to our situation of the
flag domain $D = G_0(z) \cong G_0/L_0$ where $L_0 = G_0 \cap Q_z$ is the
isotropy subgroup of $G_0$ at $z$.   Here $G_0$ replaces $A_0$,
$L_0$ replaces $B_0$\,, $E = E_\chi$ is the fiber of ${\mathbb E} \to D$ 
over $z$, the representation $\chi$ is the action of $L_0$ on $E$, and 
${\mathfrak q}_z$  replaces $\mathfrak p$.  The homogeneous holomorphic 
vector bundle structure on ${\mathbb E} \to D$ comes from an extension
$\lambda$ of $\chi$ from $L_0$ to ${\mathfrak q}_z$.   Here $\lambda$ 
integrates to $L_0$ by construction and then to $Q_z$ because $G_0 \subset G$.
Thus we have a holomorphic representation $\widetilde{\chi}$ of $Q_z$ on $E$
that extends $\chi$.  That defines the $G$--homogeneous holomorphic
vector bundle $\widetilde{\mathbb E} \to G/Q_z = Z$, and 
${\mathbb E} = \widetilde{\mathbb E}|_D$ by construction.  Statement (1) is
proved.
\s

By $G$--homogeneity of $\widetilde{\mathbb E} \to Z$ all the
$\widetilde{\mathbb E}|_C \to C$ are holomorphically equivalent. Now (2) 
follows from the construction of Leray derived sheaves.  The same 
considerations prove (3).
\end{proof}

\section {The Method of Schubert Slices} \label{schubert}
\setcounter{equation}{0}

Schubert slices play a major role in the classification part of the
above theorem (see \cite {FH} and \cite {HW}).  Here we present
the mini--version of this theory which is all that is needed for
our applications in $\S \ref {main results}$ (see \cite {H}).

\m
Let us begin by recalling that a Borel subgroup $B$ of $G$ has
only finitely many orbits in $Z$. Such an orbit ${\cal O}$
is called a Schubert cell, ${\cal O}\cong \mathbb C^{m({\cal O})}$,
and its closure $S=\bar {\cal O}={\cal O}\dot \cup Y$ is referred
to as the associated Schubert variety.  The set
${\cal S}=\{ S\}$ of all $B$--Schubert varieties freely generates
the integral homology $H_*(Z)$.

\m
If $G_0=K_0A_0N_0$ is an Iwasawa decomposition and $B\supset A_0N_0$,
then we refer to $B$ as an Iwasawa--Borel subgroup.  Recalling that
a given Borel subgroup has a unique fixed point in $Z=G/Q$, it 
follows that the Iwasawa--Borel groups are exactly those which fix a point
of the closed $G_0$--orbit in $Z$.

\m
We now prove several elementary propositions.

\begin {prop} \label {intersection property}
Let $D$ be an open $G_0$--orbit in $Z$, $C_0$ the base cycle in $D$
and $z\in D$.  Then $A_0N_0.z\cap C_0\ne \emptyset$.  In particular,
if $B$ is an Iwasawa--Borel subgroup, then every $B$--orbit ${\cal O}$
which has non--empty intersection with $D$ satisfies 
${\cal O}\cap C_0\ne \emptyset $.
\end {prop}

\begin {proof}
Since $G_0=A_0N_0K_0$ and $C_0$ is a $K_0$--orbit in $D$, it is
immediate that $D=A_0N_0.C_0$.
\end {proof}

The following holds for similar reasons.

\begin {lem} \label {tangents add up}
If $z\in D$ then $T_z(K_0.z)+T_z(A_0N_0.z)=T_zD$.
\end {lem}

\begin {prop}
If $S$ is a Schubert variety of an Iwasawa--Borel subgroup $B$ with
{\rm codim}$_\mathbb C S>q$, then $S\cap D=\emptyset $.
\end {prop}

\begin {proof}
Let ${\cal O}$ be the open dense $B$--orbit in $S$.  If 
$S\cap D\ne \emptyset $, then ${\cal O}\cap D\ne \emptyset $ as well.
Thus, by Proposition \ref {intersection property}, 
${\cal O}\cap C_0\ne \emptyset $.  But if $z\in {\cal O}\cap C_0$ then
$ \dim_\mathbb R A_0N_.z<\text{ codim}_\mathbb R K_0.z$
would contradict Lemma \ref {tangents add up}.
\end {proof}

We now come to a basic fact.

\begin {prop}
If {\em codim}$_\mathbb CS=q$ and $S\cap D\ne \emptyset $, then 
$S\cap C_0=\{z_1,\ldots z_d\} $ is non--empty, finite
and contained in the $B$--orbit ${\cal O}$.  This intersection
is transversal in the sense that

\begin {equation*}
T_{z_i}C_0\oplus T_{z_i}{\cal O}=T_{z_i}Z
\end {equation*}

for all $i$.  Furthermore, $\Sigma _i:=A_0N_0.z_i$ is open in 
${\cal O}$ and closed in $D$.
\end {prop}

\begin {proof}
Proposition \ref{intersection property} tells us that
$S\cap C_0$ is non--empty.  If $S\cap C_0$ is infinite, then it
has positive dimension at one or more of its points, contrary to 
Lemma \ref {tangents add up}.  The same argument
proves the transversality and the fact that $A_0N_0.z_i$ is open
in ${\cal O}$.  If $\Sigma _i$ were not closed in ${\cal O}$,
then we would find an $A_0N_0$--orbit of smaller dimension on its 
boundary.  By Proposition \ref{intersection property} this 
$A_0N_0$--orbit would meet $C_0$\,, contrary to Lemma \ref {tangents add up}.
\end {proof}

An $A_0N_0$--orbit $\Sigma $ as above is called a {\it Schubert slice}.
Since $[C_0]$ is non--zero in $H_*(Z;\mathbb Z)$ and
${\cal S}=\{S\}$ generates this homology, every Iwasawa--Borel
subgroup $B$ gives us $q$--codimensional Schubert varieties
with non--empty intersection $S\cap D$.  In particular, there
exist such Schubert slices.

\m 
It is known that $D$ is retractable to $C_0$ (\cite {W1}, see also
\cite {HW}) and therefore $\pi _1(D)=1$.  Let us translate this
into a statement on the isotropy groups along the base cycle.

\begin {lem} \label {isotropy splitting}
If $\Sigma $ is a Schubert slice and $z\in \Sigma \cap C_0$, then

\begin  {equation*}
\alpha :(K_0)_z\times (A_0N_0)_z\to (G_0)_z,\ (k_0,a_0n_0)\mapsto k_0a_0n_0,
\end {equation*}

is bijective.
\end {lem}

\begin {proof}
Since $K_0$ is compact and $A_0N_0$ is closed in $G_0$, it follows that
$\alpha $ is a diffeomorphism onto a closed submanifold of $G_0$.
A dimension count shows that it is open in $G_z$.  But
$G_z$ is connected, because $\pi _1(D)=1$. Thus the image of $\alpha$ is
all of $G_z$\,.
\end {proof}

We now come to the main result of this section.

\begin {thm} \label {slice}
If $\Sigma $ is a Schubert slice, then for every $C\in {\cal M}_D$
the intersection $\Sigma \cap C$ is transversal and consists of
exactly one point.
\end {thm}

\begin {proof}
First, we prove this for $C=C_0$.  Let $z_0\in C_0\cap \Sigma $.
If $z_1=a_0n_0.z_0$ were any other point in this intersection,
then $z_1=k_0^{-1}.z_0$ for some $k_0\in K_0$.  But then
$g_0=k_0a_0n_0\in (G_0)_{z_0}$ and consequently by Lemma \ref
{isotropy splitting} both $k_0$ and $a_0n_0$ fix $z_0$,
i.e., $z_1=z_0$ is the unique point in $C_0\cap \Sigma $.

\m
Now let $C\in {\cal M}_D$ be an arbitrary element of the cycle space.
Denote by ${\cal O}$ the open $B$-orbit in the Schubert variety
$S={\cal O}\dot \cup Y$ which contains $\Sigma $. Since
${\cal O}\cong \mathbb C^{m({\cal O})}$ is Stein, if $C\cap \Sigma $
were positive--dimensional then $C\cap Y\ne \emptyset $. But
codim$_\mathbb CY>q$ and therefore $Y\cap D=\emptyset $. Thus this
would be contrary to $C\subset D$.  Now $C\cap D$ is finite.

\m
If $C\cap \Sigma $ is empty we obtain a contradiction as follows.  
Let $C_t$ be a curve from $C_0$ to $C$ in ${\cal M}_D$\,.  Define 

\begin {equation*}
t_0:= \sup\{t:C_s\cap \Sigma \ne \emptyset \ \text{for all} \ s<t\}.
\end {equation*}
Then there is a sequence $\{z_n\} \subset D$ with $z_n \in C_{t_n}$ 
corresponding to $\{t_n\}$ such that  $\{t_n\} \to t_0$ and
$\{z_n\} \to z_0       \in \text{bd}(\Sigma )\subset \text{bd}(D)$. Therefore
$C_{t_0}\cap \text{bd}(D)\neq \emptyset$, contrary to $C_{t_0}\in {\cal M}_D$.

\m
This argument holds for each of the $\Sigma _i$\,.   In particular,
$C\cap S$ contains at least $d$ distinct (isolated) points.
Since $[C].[S]=d$, it follows that $\vert C\cap S\vert =d$
and that the intersection at each of these points is transversal.
\end {proof}

\section {Canonical Fibrations and DFT--Injectivity}\label {main results}
\setcounter{equation}{0}

Let $\Sigma $ be a Schubert slice defined by an Iwasawa--Borel
subgroup $B$, denote $\{z_0\} = \Sigma \cap C_0$\,,
In the context of the double fibration (\ref{flag_doublefibration}),
the projection $\nu$ carries the fiber $F = \mu ^{-1}(z_0)
=\{ (z_0,C):z_0\in C\}$ biholomorphically onto 
the analytic subset
$\{ C\in {\cal M}_D:z_0\in C\}$ of ${\cal M}_D$\,.

\m
Now consider an arbitrary element $C\in {\cal M}_Z $ with $z_0\in C$.
By definition $C=g(C_0)$ for some $g\in G$. Since $z_0\in C$, by 
adjusting $g$ by an appropriate element of $K_0$ we may assume that
$g\in Q=G_{z_0}$.  Thus
$\widetilde{F} \cong \{ C\in {\cal M}_Z:z_0\in C\} =Q.C_0$\,.  
$\widetilde{F}$ is closed in ${\cal M}_Z$\,, for if a net $\{C_i\}$
in $\widetilde{F}$ converges to $C \in {\cal M}_Z$ then $z_0 \in C$
because $z_0 \in C_i$ for each $i$, so $C \in \widetilde{F}$.
\m

The $Z$--analog of the double fibration (\ref{flag_doublefibration}) is 
given by (\ref{global_doublefibration}).
There the $Q$--orbit $\widetilde{F} \subset {\cal M}_Z $ 
is identified with the fiber $\widetilde{\mu}^{-1}(z_0)$.   In particular, 
its open subset $F = \widetilde{F} \cap \mathfrak {\cal M}_D $ is a 
closed complex submanifold on ${\cal M}_D$.

\m
By Theorem \ref {slice}, an $A_0N_0$--equivariant map 
$\varphi :{\cal M}_D\to \Sigma $ is defined by mapping $C$ to its
point of intersection with $\Sigma $.  The fiber over $z_0\in \Sigma $
is of course $F$.

\m
Let $J_0:=(A_0N_0)_{z_0}$ be the $A_0N_0$--isotropy at the base point
and note that

\begin {equation*}
A_0N_0\times _{J_0}F \to {\cal M}_D,\ \text{ defined by }
	[(a_0n_0,C)]\mapsto a_0n_0(C),
\end {equation*}

is well--defined, smooth and bijective.  Thus 
$\varphi :{\cal M}_D\to \Sigma $ is naturally identified with the smooth
$A_0N_0$--equivariant bundle 

\begin {equation*}
\pi _\Sigma :A_0N_0\times _{J_0}F \to A_0N_0/J_0=\Sigma .
\end {equation*}

In this sense, every Schubert slice defines a {\it Schubert fibration}
of the cycle space ${\cal M}_D$.

\begin {thm} \label {Schubert slice}
Let $B$ be an Iwasawa--Borel subgroup of $G$ and $\Sigma $ an
associated Schubert slice for the open orbit $D$.  Then the
fibration $\pi _\Sigma :{\cal M}_D\to \Sigma $ is a holomorphic
map on to a contractible base $\Sigma $ and diffeomorphically
realizes ${\cal M}_D$ as the product $\Sigma \times F$.
\end {thm}

\begin {proof}  Let $J:= B_{z_0}$\,.
The inclusions $A_0N_0 \hookrightarrow B$ and 
$F \hookrightarrow \widetilde{F}$ together define a map
$A_0N_0\times_{J_0}  F \hookrightarrow B\times _J F$.  That map realizes 
$A_0N_0\times _{J_0}F \cong {\cal M}_D$
as an open subset of $B\times _J\widetilde{F}$.  The latter is fibered over the
open $A_0N_0$--orbit $\Sigma $ in ${\cal O}=B.z_0$ by the natural holomorphic
projection $\pi :B\times _J\widetilde{F} \to B/J$.  Since $\pi _\Sigma $ is the
restriction $\pi \vert_{{\cal M}_D}$\,, it follows that $\pi _\Sigma $
is holomorphic as well.

\s
The fact that $\Sigma $ is a cell follows form the simple connectivity
of the solvable group $A_0N_0$ and the fact that it is acting algebraically.
\end {proof}

Recall the notation: $\Omega_\mu^r(\mathbb E) \to \mathfrak X_D$ is the
sheaf of relative $\mu^*\mathbb E$--valued holomorphic $r$--forms on 
$\mathfrak X$ with respect to $\mu : \mathfrak X_D \to D$.

\begin {cor} \label {DFT--injectivity}
Suppose that $\mathbb E\to D$ a holomorphic $G_0$--homogeneous vector bundle
which is sufficiently negative so that $H^p(C; \Omega_\mu^r(\mathbb E)|_C) = 0$ 
for $p < q$, and $r \geqq 1$.  Then the double fibration transform

\begin {equation*}
{\mathcal P}:H^q(D,{\cal O}(\mathbb E)) \to H^0(M,{\cal O}(\mathbb E'))
\end {equation*}

is injective.
\end {cor}

\begin {proof}
${\cal M}_D$ is contractible by Proposition \ref {contractible}.  Since 
$\Sigma $ is likewise contractible and ${\cal M}_D$ is diffeomorphic
to $\Sigma \times F$, it follows that $F$ is cohomologically trivial.
The Buchdahl conditions (\ref{buchdahl_conditions}) follow.
Proposition \ref {buchdahl_theorem} now says that 
(\ref{pull1}) is an isomorphism for all $r$.  Composing with coefficient
morphisms, the maps (\ref{pull2}) also are isomorphisms.
By Theorem \ref{cycle spaces} we know that the conditions 
(\ref{proper_stein}) are
satisfied.  The assertion now follows from Theorem \ref{general_injective}.
\end {proof}
\m

With a bit more work one can see that the fiber $F$ of ${\cal M}_D \to \Sigma$
is contractible, not just cohomologically trivial.  We thank Peter Michor for
showing us the following result for the $C^\infty$ category, from which
contractibility of $F$ is immediate.  His argument is based on the existence
of a complete Ehresmann connections for smooth fiber bundles.
 
\begin{prop} \label{fiber_contractible}
Let $p: M \to S$ be a smooth fiber bundle with fiber $F = p^{-1}(s_0)$.
If both $M$ and $S$ are contractible then $F$ is contractible.
\end{prop}
 
\begin {proof} Since $S$ is contractible and smooth, approximation
gives us a smooth contraction $h : [0,1] \times S \to S$; here $h(0,s) = s$
and $h(1,s) = s_0$\,.  Following \cite[\S 9.9]{KMS} the bundle $p: M \to S$ has
a complete Ehresmann connection.  Completeness means that every smooth curve
in $S$ has horizontal lifts to $M$.  If $m \in M$ let $t \mapsto H(t,m)$
denote the horizontal lift of $t \mapsto h(t,p(m))$ such that $H(0,m) = m$.
Note $H(1,m) \in F$.  Fix a base point $m_0 \in M$ and a smooth contraction
$I : [0,1] \times M \to M$ of $M$ to $m_0$\,; if $m \in M$ then $I(0,m) = m$
and $I(1,m) = m_0$.  Denote $f_0 = H(1,m_0) \in F$.  Define
$J : [0,1] \times F$ by $J(t,f) = H(t,I(t,f))$, so $J(0,f) = f$ and
$J(1,f) = f_0$\,.  Thus $J$ is a contraction of $F$ to $f_0$\,, and so
$F$ is contractible.
\end{proof}

\vskip 1 cm

\centerline{\begin{tabular}{ll}
ATH: & JAW: \\
Fakult\" at f\" ur Mathematik & Department of Mathematics \\
Ruhr--Universit\" at Bochum & University of California \\
D-44780 Bochum, Germany & Berkeley, California 94720--3840, U.S.A. \\
                               &                                      \\
{\tt ahuck@cplx.ruhr-uni-bochum.de} & {\tt jawolf@math.berkeley.edu}
\end{tabular}}

\end {document}